# On nilpotency in pre - Lie algebras

**L.A. Kurdachenko, I Ya. Subbotin**

**Abstract:** The authors consider several variants of the notion of nilpotency in pre-Lie algebras and investigate the properties of nilpotent ideals.

A left pre - Lie algebra A over a field F is a vector space with binary operations

$$(x, y) \longrightarrow x + y, (x, y) \longrightarrow x \cdot y, x, y \in A,$$

which satisfies the following conditions

$$(x \cdot y) \cdot z - x \cdot (y \cdot z) = (y \cdot x) \cdot z - y \cdot (x \cdot z),$$
$$(x + y) \cdot z = x \cdot z + y \cdot z, x \cdot (y + z) = x \cdot y + x \cdot z$$

for all elements $x, y, z \in A$.

Define the associator of elements x, y, z by the following formula

$$(x, y, z) = (x \cdot y) \cdot z - x \cdot (y \cdot z).$$

Then the first equality looks as $(x, y, z) = (y, x, z)$. In other words, in the pre-Lie algebras the associator is symmetric in the left two variables x, y. For this reason, in E. Vinberg's paper [VE1963], these algebras were called left symmetric algebras. In that paper, the author studies convex homogeneous cones. This concept has been applied in numerous studies related to geometry. The term a pre – Lie algebra first appeared in the paper by M. Gerstenhaber [GM1963]. Pre-Lie algebras, under other names, can also be found in the papers by B. Kuperschmidt [KB1994] and J.-L. Koszul [KJL1961]. Since then, pre-Lie algebras have become an active area of investigation due to their deep connections with affine geometry, deformation theory, operads, Lie theory, combinatorics, and mathematical physics.

By symmetry, a **right pre - Lie algebra** A over a field F is a vector space with binary operations

$$(x, y) \longrightarrow x + y, (x, y) \longrightarrow x \cdot y, x, y \in A,$$

such that (A, +) is an abelian group satisfying the following conditions

$$(x \cdot y) \cdot z - x \cdot (y \cdot z) = (x \cdot z) \cdot y - x \cdot (z \cdot y),$$
$$(x + y) \cdot z = x \cdot z + y \cdot z, x \cdot (y + z) = x \cdot y + x \cdot z$$

for all elements $x, y, z \in A$.

Let L be a right pre–Lie algebra. Define the new operation $\bullet$ by the rule

$$x \bullet y = y \cdot x \text{ for all } x, y \in A.$$

It is possible to show that L is a left pre – Lie algebra by the operations + and $\bullet$. Thus, it shows that we can consider only left pre-Lie algebras.

The importance of pre-Lie algebras is further illustrated by the fact that many classical algebraic structures can be naturally related to them. Every associative algebra is a pre-Lie algebra. Furthermore, important subclasses such as assosymmetric algebras and Novikov algebras arise by imposing additional identities.

An algebra L is called **bi – symmetric algebra or assosymmetric algebra** if A is left and right pre – Lie algebra.

Such algebras were introduced by Kleinfeld in [KZ1957] under the name *assosymmetric algebras*, preceding the systematic study of pre-Lie algebras.

A pre - Lie algebra L over a field F is called a **Novikov algebra** if it satisfies the following condition

$$(x \bullet y) \bullet z = (x \bullet z) \bullet y$$

for all elements $x, y, z \in L$.

Novikov algebras arise in various areas of mathematics and mathematical physics. In particular, they were introduced in the study of Hamiltonian operators related to the integrability of certain nonlinear partial differential equations [GD1980]. They also appear naturally in the theory of Poisson brackets of hydrodynamic type [BN1985]. One of the first fundamental results on the structure of Novikov algebras was obtained by Zelmanov [ZE1987].

Structural properties of Novikov algebras have consequently received considerable attention and continue to motivate the development of the general theory of pre-Lie algebras.

A pre–Lie algebras L is called **commutative** if $y \cdot x = x \cdot y$ for all elements $x, y \in L$. Note also that every commutative pre – Lie algebra is associative [KS2025].

A pre - Lie algebras L is called **abelian** if $y \cdot x = 0$ for all elements $x, y \in L$.

The study of nilpotency conditions constitutes one of the central directions of research in many areas of algebra, including group theory, associative rings and algebras, Lie algebras, Leibniz algebras, and other algebraic structures. The same theme plays an important role in the theory of pre-Lie algebras.

In the associative rings this notion looks as the following. An associative ring R is called nilpotent if there exists a positive integer k such that the product of every k elements of R is zero. In other words the subset

$$R^{ok} = \{ a_1a_2 \dots a_k \mid a_1, a_2, \dots, a_k \in R \} = \langle 0 \rangle.$$

In an associative ring, the product of $k$ elements does not depend on the placement of parentheses. The situation is considerably more complicated in nonassociative algebras. In this case, the value of a product generally depends on the placement of brackets, and several nonequivalent notions of nilpotency may arise.

Consequently, different definitions are possible, each reflecting a particular way of constructing iterated products. One of the principal objectives of the theory is to understand the relationships among these definitions and determine which properties of classical nilpotent structures remain valid in the nonassociative setting.

For pre-Lie algebras, this phenomenon leads naturally to the notions of left nilpotency, right nilpotency, and central nilpotency. These concepts are defined through different series associated with the algebra and need not coincide a priori. The resulting structure theory is therefore richer than in the associative case and presents a number of difficulties absent from classical theories.

Let $L$ be a pre – Lie algebra. Put $L^{(1)} = L$, recursively $L^{(\alpha + 1)} = L^{(\alpha)} \cdot L$ for all ordinals $\alpha$, and $L^{(\lambda)} = \cap_{\mu < \lambda} L^{(\mu)}$ for limit ordinals $\lambda$. **Proposition 1.2** of paper [KS2025] shows that $L^{(\alpha)}$ is an ideal of $L$.

We say that a pre-Lie algebra $L$ is a **right nilpotent** if $L^{(n)} = \langle 0\rangle$ for some positive integer $n$.

Similarly, put $L^1 = L$, recursively $L^{\alpha + 1} = L \cdot L^{\alpha}$ for all ordinals $\alpha$, and $L^{\lambda} = \cap_{\mu < \lambda} L^{\mu}$ for limit ordinals $\lambda$. In general, $L$ is not an ideal, **Proposition 1.2 in** [KS2025] **shows merely** that $L^{\alpha}$ is a left ideal of $L$.

We say that a pre-Lie algebra $L$ is a **left nilpotent** if $L^n = \langle 0\rangle$ for some positive integer $n$.

Define the **lower central series** of a pre – Lie algebra $L$

$$L = \gamma_1(L) \geq \gamma_2(L) \geq \ldots \geq \gamma_\alpha(L) \geq \gamma_{\alpha+1}(L) \geq \ldots \gamma_\delta(L)$$

by the following rule: $\gamma_1(L) = L$, $\gamma_2(L) = L^{(2)}$, $\gamma_3(L) = L^{(2)} \cdot L + L \cdot L^{(2)} = \gamma_2(L) \cdot L + L \cdot \gamma_2(L)$, recursively $\gamma_{\alpha+1}(L) = \gamma_\alpha(L) \cdot L + L \cdot \gamma_\alpha(L)$ for all ordinals $\alpha$, and $\gamma_\lambda(L) = \cap_{\mu < \lambda} \gamma_\mu(L)$ for limit ordinals $\lambda$. The last term $\gamma_\delta(L)$ is called the **lower hypocenter of L.** We have $\gamma_\delta(L) = \gamma_\delta(L) \cdot L + L \cdot \gamma_\delta(L)$.

Using again **Proposition 1.2** of a paper [KS2025] we obtain that $\gamma_\alpha(L)$ is an ideal of $L$.

Denote by $\mathbf{lzl}(L)$ the length of lower central series of $L$.

Let $L$ be a pre-Lie algebra over a field $F$, $M$ be a non-empty subset of $L$ and $H$ be a subalgebra of $L$. Put

$$\mathbf{Ann}_H^{\mathbf{left}}(M) = \{ a \in H \mid a \cdot u = 0 \text{ for each element } u \in M \},$$
$$\mathbf{Ann}_H^{\mathbf{right}}(M) = \{ a \in H \mid u \cdot a = 0 \text{ for each element } u \in M \}.$$

The subset $\mathbf{Ann}_H^{\mathbf{left}}(M)$ is called the **left annihilator of M in subalgebra H**; the subset $\mathbf{Ann}_H^{\mathbf{right}}(M)$ is called the **right annihilator of M in a subalgebra H.** The intersection

$$\mathbf{Ann}_H(M) = \{ a \in H \mid a \cdot u = 0 = u \cdot a \text{ for each element } u \in M \}$$

is called the **annihilator of M in subalgebra H.** It is clear that $\mathbf{Ann}_H(M) = \mathbf{Ann}_H^{\mathbf{left}}(M) \cap \mathbf{Ann}_H^{\mathbf{right}}(M)$.

Let $L$ be a pre-Lie algebra, put

$$\zeta_{\mathbf{ann}}^{\mathbf{left}}(L) = \{ x \in L \mid x \cdot y = 0 \text{ for each element } y \in L \} = \mathbf{Ann}_H^{\mathbf{left}}(A),$$
$$\zeta_{\mathbf{ann}}^{\mathbf{right}}(L) = \{ x \in L \mid y \cdot x = 0 \text{ for each element } y \in L \}) = \mathbf{Ann}_H^{\mathbf{right}}(A)$$

and

$$\zeta_{\mathbf{ann}}(L) = \{ x \in L \mid x \cdot y = y \cdot x = 0 \text{ for each element } y \in L \}) = \mathbf{Ann}_H(A).$$

The subset $\zeta^{left}(L)$ is called the **left annihilation center of L**, the subset $\zeta^{right}(M)$ is called the **right annihilation center of L**, and the subset $\zeta(L)$ is called **the annihilation center of L.**

**Proposition 1.5** in [KS2025] shows that $\zeta^{left}(L)$ and $\zeta(L)$ are ideals of L.

Define the **upper central series**

$$\langle 0\rangle = \zeta_0(L) \leq \zeta_1(L) \leq \zeta_2(L) \leq \ldots \leq \zeta_\alpha(L) \leq \zeta_{\alpha+1}(L) \leq \ldots \zeta_\gamma(L) = \zeta_\infty(L)$$

of a pre-Lie algebra L by the following rule: $\zeta_1(L) = \zeta_{ann}(L)$ is the annihilation center of L, recursively $\zeta_{\alpha+1}(L)/\zeta_\alpha(L) = \zeta_{ann}(L/\zeta_\alpha(L))$ for all ordinals $\alpha$, and $\zeta_\lambda(L) = \cup_{\mu<\lambda} \zeta_\mu(L)$ for limit ordinals $\lambda$. Each term of this series is an ideal of L. The last term $\zeta_\infty(L)$ of this series is called the **upper hypercenter of L.** Denote by **uzl**(L) the length of the upper central series of L.

If $L = \zeta_\infty(L)$, then L is said to be **a hypercentral pre-Lie algebra.**

Let L be a pre-Lie algebra over a field F, A, B be ideals of L such that $B \leq A$. Then the factor A/B is called **central in L** if $A/B \leq \zeta(L/B)$. This means that $A \cdot L$ and $L \cdot A \leq B$, and hence, $A \cdot L + L \cdot A \leq B$.

We note that every factor of the lower and upper central series of L is central in L.

**Proposition 1.7** of paper [KS2025] shows that if a pre-Lie algebra L has a finite series of ideals, whose factors are central in L, then the lower and upper central series of L are finite and **uzl**(L) = **lzl**(L). Therefore, as usual, we shall call a pre-Lie algebra L **nilpotent** if it possesses a finite series of ideals whose successive factors are central in L.

One of the first important results concerning nilpotent ideals is the assertion that the sum of two nilpotent ideals is again a nilpotent ideal. An analogous statement holds for groups, where it is known as Fitting's theorem. This result constitutes the first step toward the construction of radicals.

For this reason, it is natural to ask whether analogous results remain valid for pre-Lie algebras. The answer is far from obvious. The coexistence of several notions of nilpotency, together with the asymmetry of the left and right multiplication operators, suggests that different phenomena may arise depending on the chosen definition. Consequently, understanding the behavior of sums of nilpotent ideals becomes a central problem.

The primary purpose of the present paper is to investigate this question. We study several notions of nilpotency in pre-Lie algebras and analyze the structure of sums of nilpotent ideals. Our results demonstrate that the situation differs substantially according to whether left or right nilpotency is considered. On the one hand, a number of positive analogues of classical Fitting-type theorems can be obtained. On the other hand, right nilpotency exhibits additional complications that require more delicate methods.

We now summarize the main results of the paper.

Our first main result concerns left nilpotent ideals. It provides an initial step toward a general Fitting-type theorem for left nilpotent ideals. The proof reveals the mechanisms by which left nilpotency is inherited by the sum of nilpotent ideals and suggests possible extensions to broader classes of nilpotent ideals.

**THEOREM A.** *Let $L$ be a pre – Lie algebra, $B$, $C$ be ideals of $L$ such that $L = B + C$. Suppose that $B^k = C^3 = <0>$. If $B$ is abelian, then $L^4 = <0>$. If $k > 2$, then $L^{2k} = <0>$.*

The proof of this theorem reveals the mechanisms through which left nilpotency propagates from the summands to their sum. It also indicates a possible approach to proving the more general statement that the sum of two left nilpotent ideals is again left nilpotent, although it is clear that the necessary computations would be substantially more involved.

The next theorem shows that the corresponding problem for right nilpotency is considerably more subtle.

**THEOREM B.** *Let $L$ be a pre – Lie algebra, $B$, $C$ be ideals of $L$ such that $L = B + C$. If $B^{(3)} = C^{(3)} = <0>$, then $B \cap C$ includes an abelian ideal $J$ such that $(L/J)^{(6)} = <0>$.*

Thus, in contrast to the left nilpotent case, one cannot expect a direct analogue of Fitting's theorem. The theorem nevertheless shows that the obstruction is controlled by a suitable abelian ideal.

The final theorem of the paper provides detailed information on the sum of two nilpotent ideals of central length at most three.

**THEOREM C.** *Let $L$ be a pre – Lie algebra, $B$, $C$ be ideals of $L$ such that $L = B + C$, $B = \zeta_2(B)$, $C = \zeta_2(C)$. Then the following assertions hold*

*(i)* $\zeta_1(B) \cap \zeta_1(C) \leq \zeta_1(L)$.

*(ii)* $B \bullet (C \bullet D) \leq \zeta_1(B) \cap \zeta_1(C)$, $C \bullet (B \bullet D) \leq \zeta_1(B) \cap \zeta_1(C)$.

*(iii)* $L \bullet D = B \bullet D + C \bullet D \leq (D \cap \zeta_1(B)) + (D \cap \zeta_1(C)) \leq \zeta_1(D)$,

*(iv)* $L \bullet (B \bullet D + C \bullet D) \leq B \bullet (C \bullet D) + C \bullet (B \bullet D) \leq \zeta_1(B) \cap \zeta_1(C) \leq \zeta_1(L)$, *in particular,* $B \bullet (C \bullet D) + C \bullet (B \bullet D)$ *is an ideal of* $L$.

*(v)* $L \bullet (L \bullet (L \bullet D) = <0>$.

*(vi)* $L^6 = <0>$.

*(vii)* $D \bullet L \leq D \bullet B + D \bullet C \leq (D \cap \zeta_1(B)) + (D \cap \zeta_1(C)) \leq \zeta_1(D)$,

*(viii)* $(D \cap \zeta_1(B)) + (D \cap \zeta_1(C))$ *is an abelian ideal of* $L$.

*(ix)* $\zeta_1(D)$ *is an ideal of* $L$.

*(x)* $(D \bullet L) \bullet L \leq (D \bullet C) \bullet B + (D \bullet B) \bullet C \leq (D \cap \zeta_1(B)) + (D \cap \zeta_1(C)) \leq \zeta_1(D)$, $((D \bullet L) \bullet L) \bullet L \leq (D \bullet C) \bullet B + (D \bullet B) \bullet C \leq \zeta_1(D)$.

*(xi)* $L^{(5)} \leq (D \cap \zeta_1(B)) + (D \cap \zeta_1(C)) \leq \zeta_1(D)$.

*(xii) The factor – algebra* $L/((D \cap \zeta_1(B)) + (D \cap \zeta_1(C)))$ *is nilpotent and* **uzl**$(A/((D \cap \zeta_1(\ast, B)) + (D \cap \zeta_1(\ast, C)))) \leq 3$.

Thus, in contrast to the left nilpotent case, one cannot expect a direct analogue of Fitting's theorem. The theorem nevertheless shows that the obstruction is governed by a suitable abelian ideal.

Our final theorem provides a detailed description of the sum of two nilpotent ideals whose central length does not exceed three.

## Preliminaries

If $L$ is a pre – Lie algebra and $K$ be an ideal of $L$, then $K \bullet L$ is an ideal of $L$ [KS2025, Proposition 1.2]. However, in general, $L \bullet K$ is a left ideal of $L$. But

**1.1. PROPOSITION.** *Let $L$ be a pre – Lie algebra and $K$ be an ideal of $A$. If $\zeta^{left}(L)$ includes $K$, then $L \bullet K$ is an ideal of $L$.*

**PROOF.** By **Proposition 1.2** of [KS2025] $L \bullet J$ is a left ideal of $L$.

Let $a$ be an arbitrary element of $L \bullet K$, then $a = \sum_{1 \leq j \leq n} x_j \cdot b_j$, where $x_j \in L$, $b_j \in K$, $1 \leq j \leq n$. If $y$ is an arbitrary element of $L$, then

$$a \cdot y = (\sum_{1 \leq j \leq n} x_j \cdot b_j) \cdot y = \sum_{1 \leq j \leq n} ((x_j \cdot b_j) \cdot y).$$

We have $(x_j \cdot b_j) \cdot y - x_j \cdot (b_j \cdot y) = (x_j \cdot y) \cdot b_j - x_j \cdot (y \cdot b_j)$. Inclusion $K \leq \zeta^{left}(L)$ implies that $b_j \cdot y = 0$, and hence $x_j \cdot (b_j \cdot y) = 0$. Since $K$ is an ideal, $y \cdot b_j \in K$, we have $x_j \cdot (y \cdot b_j) \in L \bullet K$. Since $(x_j \cdot y) \cdot b_j \in L \bullet K$, we obtain that $(x_j \cdot b_j) \cdot y \in L \bullet K$. It is true for each $j$, $1 \leq j \leq n$, and therefore $a \cdot y \in L \bullet K$.

**1.2. PROPOSITION.** *Let $L$ be a pre – Lie algebra and $S = \zeta_{ann}{}^{left}(L)$. If $L^{n+1} = <0>$, then $S$ has a finite series of ideals of $L$ whose length is at most $n$ and whose factors are central in $L$.*

**PROOF.** Let

$$L = L^1 \geq L^2 \geq \ . \ . \ . \geq L^{n+1} = <0>$$

be the lower left central series of $L$. Put $S_2 = L \bullet S$. By **Proposition 1.1** $S_2$ is an ideal of $L$. Since $S$ is an ideal, $L \bullet S \leq S$. Clearly $S_2 \leq L^2$ and hence $S_2 \leq S \cap L^2$. We have $L \bullet S \leq S_2$ and $<0> = S \bullet L$, so that the factor $S/S_2$ is central in $L$.

Put $S_3 = L \bullet S_2$. By **Proposition 1.1** $S_3$ is an ideal of $L$. The inclusion $S_2 \leq L^2$ implies that $S_3 = L \bullet S_2 \leq L \bullet L^2 = L^3$, so that $L_3 \leq A^3 \cap S$. We have $L \bullet S_2 \leq S_3$ and $<0> = S_2 \bullet L_2$, so that the factor $S_2/S_3$ is central in $L$.

In general, put $S_{j+1} = L \bullet S_j$ for every positive integer $j$. Using similar arguments and ordinary induction we obtain that $S_j$ is an ideal of $A$, $S_j \leq L^j \cap S$ and the the factor $S_j/S_{j+1}$ for every positive integer $j$ is central in $L$. The equality $<0> = A^{n+1}$ implies that $S_{n+1} = <0>$.

**1.3. COROLLARY.** *Let $L$ be a pre – Lie algebra and let $n$, $k$ be positive integers such that $L^{(n)} = <0> = L^k$. Then $L$ is nilpotent and and $\mathbf{lzl}(L) \leq nk$.*

**PROOF.** Since $L^{(n)} = <0>$, $L$ has a finite series of ideals

$$L = L^{(1)} \geq L^{(2)} \geq \ . \ . \ . \geq L^{(n-1)} \geq L^{(n)} = <0>.$$

Then $L^{(j)}/L^{(j+1)} \leq \zeta_{ann}{}^{left}(L/L^{(j+1)})$, $1 \leq j \leq n - 1$. ${}_{ann}{}^{left}(L)$. By **Proposition 1.2** $L^{(j)}/L^{(j+1)}$ has a finite series of ideals whose length is at most $k$ and whose factors are central in $L$. Then $L$ has has a finite series of ideals whose length is at most $nk$ and whose factors are central in $L$. Then **Proposition 1.7** of a paper [KS2025] implies that $L$ is nilpotent and and $\mathbf{lzl}(L) \leq nk$.

**1.4. LEMMA.** *Let $L$ be a pre – Lie algebra, $B, C$ are ideals of $L$ such that $L = B + C$. If $B \cap C = <0>$, then $L^n = B^n + C^n$, $L^{(n)} = B^{(n)} + C^{(n)}$, $\zeta_n(L) = \zeta_n(B) + \zeta_n(C)$ for every positive integer $n$.*

**PROOF.** If $w$ be arbitrary element of $L^2$, then $w = \sum_{1 \leq j \leq k} x_j \cdot y_j$, where $x_j, y_j \in L$, $1 \leq j \leq k$. In turn $x_j = b_{1j} + c_{1j}$, $y_j = b_{2j} + c_{2j}$ for some elements $b_{1j}, b_{2j} \in B$, $c_{1j}, c_{2j} \in C$, $1 \leq j \leq k$. We have

$$x_j \cdot y_j = (b_{1j} + c_{1j}) \cdot (b_{2j} + c_{2j}) = b_{1j} \cdot b_{2j} + b_{1j} \cdot c_{2j} + c_{1j} \cdot b_{2j} + c_{1j} \cdot c_{2j}.$$

Since $B, C$ are the ideals, $b_{1j} \cdot c_{2j}, c_{1j} \cdot b_{2j} \in B \cap C = <0>$, so that $b_{1j} \cdot c_{2j} = c_{1j} \cdot b_{2j} = 0$. Thus $x_j \cdot y_j = b_{1j} \cdot b_{2j} + c_{1j} \cdot c_{2j} \in B^2 + C^2$. Then $w = \sum_{1 \leq j \leq k} x_j \cdot y_j \in B^2 + C^2$.

Suppose that we have already proved that $L^j = B^j + C^j$ for all $j < m$. Every element $v$ of $L^m$ has a form $v = \sum_{1 \leq j \leq k} u_j \cdot z_j$, where $u_j \in L$, $z_j \in L^{m-1}$, $1 \leq j \leq k$. We have $u_j = b_{3j} + c_{3j}$ for some elements $b_{3j} \in B$, $c_{3j} \in C$, $1 \leq j \leq k$. The equality $L^{m-1} = B^{m-1} + C^{m-1}$ implies that $z_j = b_{4j} + c_{4j}$ where $b_{4j} \in B^{m-1}$, $c_{4j} \in C^{m-1}$, $1 \leq j \leq k$. Now we obtain

$$v = \sum_{1 \leq j \leq k} u_j \cdot z_j = \sum_{1 \leq j \leq k} (b_{3j} + c_{3j})(b_{4j} + c_{4j}) = \sum_{1 \leq j \leq k}(b_{3j}b_{4j} + b_{3j}c_{4j} + c_{3j}b_{4j} + c_{3j}c_{4j}).$$

We have that $b_{3j}c_{4j}, c_{3j}b_{4j} \in B \cap C = <0>$, so that $b_{3j} \cdot c_{4j} = c_{3j} \cdot b_{4j} = 0$ and we obtain

$$v = \sum_{1 \leq j \leq k} u_j \cdot z_j = \sum_{1 \leq j \leq k}(b_{3j}b_{4j} + c_{3j}c_{4j}) = \sum_{1 \leq j \leq k} b_{3j}b_{4j} + \sum_{1 \leq j \leq k} c_{3j}c_{4j} \in B^m + C^m.$$

Using similar arguments we can prove an equality $L^{(n)} = B^{(n)} + C^{(n)}$.

Let $z$ be the arbitrary element of $\zeta_1(B)$, $x$ be an arbitrary element of $A$. Then $x = b + c$ where $b \in B$, $c \in C$. We have $z \cdot x = z \cdot (b + c) = z \cdot b + z \cdot c$. Since $B, C$ are ideals, $z \cdot c \in B \cap C = <0>$, so that $z \cdot x = z \cdot b = 0$. In a similar way we obtain that $x \cdot z = 0$. It follows that $z \in \zeta_1(L)$, and it implies that $\zeta_1(B) \leq \zeta_1(L)$. Similarly, $\zeta_1(C) \leq \zeta_1(L)$, and hence $\zeta_1(B) + \zeta_1(C) \leq \zeta_1(L)$. Conversely, let $y \in \zeta_1(L)$, $y = y_1 + y_2$ where $y_1 \in B$, $y_2 \in C$. If $b$ is an arbitrary element of $B$, then

$$0 = y \cdot b = (y_1 + y_2) \cdot b = y_1 \cdot b + y_2 \cdot b.$$

Since $B, C$ are ideals, $y_2 \cdot b \in B \cap C = <0>$, so that $y_1 \cdot b = 0$. In a similar way we obtain that $b \cdot y_1 = 0$, so that $y_1 \in \zeta_1(B)$. Similarly, we can prove that $y_2 \in \zeta_1(C)$. This proves the equality $\zeta_1(B) + \zeta_1(C) = \zeta_1(L)$.

Using ordinary induction and the same arguments we prove the equality $\zeta_n(L) = \zeta_n(B) + \zeta_n(C)$ for each positive integer $n$.

**1.5. COROLLARY.** *Let $L$ be a pre – Lie algebra, $B, C$ are ideals of $L$ such that $L = B + C$ and $B \cap C = <0>$.*

*If $B^n = C^k = <0>$, then $A^m = <0>$ where $m = \max \{n, k\}$,*
*If $B^{(n)} = C^{(k)} = <0>$, then $A^{(m)} = <0>$ where $m = \max \{n, k\}$,*
*If $B = \zeta_n(B)$, $C = \zeta_k(C)$, then $A = \zeta_m(A)$ where $m = \max \{n, k\}$.*

**1.6. LEMMA.** *Let $L$ be a pre – Lie algebra and $B$ is an ideal of $L$. Then $\zeta_{ann}{}^{right}(L)$ includes both products $\zeta(B) \bullet L$, $L \bullet \zeta(B)$.*

**PROOF.** Let $z$ be an arbitrary element of $\zeta(B)$, $b$ be an arbitrary element of $B$ and $x$ be an arbitrary element of $L$. We have

$$(b \cdot z) \cdot x - b \cdot (z \cdot x) = (z \cdot b) \cdot x - z \cdot (b \cdot x).$$

Since $z \in \zeta(B)$, $z \cdot b = 0 = b \cdot z$. Since B is an ideal, $b \cdot x \in B$, and hence $z \cdot (b \cdot x) = 0$. Thus, we obtain $b \cdot (z \cdot x) = 0$, so that $z \cdot x \in \zeta_{ann}{}^{right}(L)$. Similarly,

$$(b \cdot x) \cdot z - b \cdot (x \cdot z) = (x \cdot b) \cdot z - x \cdot (b \cdot z).$$

Again, we have $(b \cdot x) \cdot z = (x \cdot b) \cdot z = x \cdot (b \cdot z) = 0$, and hence $z \cdot (b \cdot x) = 0$. Thus, we obtain $b \cdot (x \cdot z) = 0$, so that $x \cdot z \in \zeta_{ann}{}^{right}(L)$.

## Sums of left and right nilpotent ideals

**2.1. PROPOSITION.** *Let L be a pre – Lie algebra, B, C are ideals of L. If B is abelian and C is nilpotent, then C + B is also nilpotent and* $\mathbf{zl}(C + B) \le 2\mathbf{zl}(C)$.

**PROOF.** If $B \cap C = \langle 0 \rangle$, then the result follows from **Corollary 1.5**.

Therefore, suppose that $B \cap C = D \ne \langle 0 \rangle$. In this case we use induction on $\mathbf{zl}(C) = k$. Let

$$\langle 0 \rangle = C_0 \le C_1 \le \dots \le C_k = C$$

be the upper central series of C, $S_j = C_j \cap B$, $0 \le j \le k$. Let d be arbitrary element of $S_1$ and x be arbitrary element of B + C, then $x = b + c$ where $b \in B, c \in C$. Since $d \in C_1$, $d \cdot c = c \cdot d = 0$. The fact that B is abelian implies that $d \cdot b = b \cdot d = 0$. Then

$$d \cdot (b + c) = d \cdot b + d \cdot c = 0, (b + c) \cdot d = b \cdot d + c \cdot d = 0.$$

It follows that $d \in \zeta(B + C)$, so that $S_1 \le \zeta_1(B + C)$. Using similar arguments and ordinary induction we obtain that $S_j \le \zeta_j(B + C)$, $0 \le j \le k$. In particular, $D = S_k \le \zeta_k(B + C)$. Applying **Corollary 1.5** to the factor – algebra L/D we obtain that L/D is nilpotent and $\mathbf{zl}(L/D) \le k$. It follows that L is nilpotent and $\mathbf{zl}(L) \le 2k$.

**2.2. COROLLARY.** *Let L be a pre – Lie algebra, B, C are abelian ideals of A. Then their sum C + B is nilpotent and* $\mathbf{zl}(C + B) \le 2$.

**2.3. PROPOSITION.** *Let L be a pre – Lie algebra, L = B + C where B, C the ideals of L.*

*If* $B^n = \langle 0 \rangle$ *and C is an abelian ideal, then* $L^{2n-1} = \langle 0 \rangle$.

*If* $B^{(n)} = \langle 0 \rangle$ *and C is an abelian ideal, then* $L^{(2n-1)} = \langle 0 \rangle$.

**PROOF.** Let $D = B \cap C$. **Corollary 1.5** implies that $(L/D)^n = \langle 0 \rangle$. The equality $(L/D)^n = (L^n + D)/D$ implies that $L^n \le D$.

Let x be an arbitrary element of A, d be an arbitrary element of D, then $x = b + c$ where $b \in B, c \in C$. We have $x \cdot d = (b + c) \cdot d = b \cdot d + c \cdot d$. Since C is abelian, $c \cdot d = 0$ and $x \cdot d = b \cdot d$. The inclusion $D \le B$ implies that $b \cdot d \in B^2$. On the other hand, the fact that D is an ideal implies that $b \cdot d \in D$, so that $b \cdot d \in B^2 \cap D$. Thus, we obtain that $L \bullet D \le B^2 \cap D$. It follows that $L^{n+1} = L \bullet L^n \le L \bullet D \le B^2 \cap D$. Using similar argument, we prove that

$$L^{n+2} \leq B^3 \cap D, \ldots, L^{n+n-1} \leq B^n \cap D = \langle 0\rangle.$$

The proof of a second assertion is similar.

If $L$ is a pre – Lie algebra, then the equality $(b \cdot c) \cdot d – b \cdot (c \cdot d) = (c \cdot b) \cdot d - c \cdot (b \cdot d)$ implies

$$b \cdot (c \cdot d) = (b \cdot c) \cdot d – (c \cdot b) \cdot d + c \cdot (b \cdot d).$$

We will use this equality frequently throughout the paper.

**2.4. LEMMA.** *Let $L$ be a pre – Lie algebra, B, C are the ideals of L, $D = B \cap C$. Then*

$$B \bullet (C \bullet D) \leq (B \bullet D) \cap (C \bullet D).$$

**PROOF.** Every element $x$ of $B \bullet (C \bullet D)$ can be represented in the form $x = \sum_{1 \leq j \leq n} b_j \cdot u_j$ for some elements $b_j \in B$, $u_j \in C \bullet D$, $1 \leq j \leq n$. In turn, element $u_j \in C \bullet D$ has a form $\sum_{1 \leq k \leq m} c_{jkj} \cdot d_{jk}$ for some elements $c_{jk} \in C$, $d_{jk} \in D$, $1 \leq j \leq m$. Thus, we have

$$x = \sum_{1 \leq j \leq n} b_j \cdot u_j = \sum_{1 \leq j \leq n} b_j \cdot (\sum_{1 \leq k \leq m} c_{jkj} \cdot d_{jk}) =$$
$$\sum_{1 \leq j \leq n} \sum_{1 \leq k \leq m} (b_j \cdot (c_{jkj} \cdot d_{jk})).$$

Consider an element $b_j \cdot (c_{jkj} \cdot d_{jk})$. The fact that $D$ is an ideal implies that $c_{jkj} \cdot d_{jk} \in D$, so that $b_j \cdot (c_{jkj} \cdot d_{jk}) \in B \bullet D$. On the other hand, the equality

$$(b_j \cdot c_{jkj}) \cdot d_{jk} – b_j \cdot (c_{jkj} \cdot d_{jk}) = (c_{jkj} \cdot b_j) \cdot d_{jk} – c_{jkj} \cdot (b_j \cdot d_{jk})$$

implies the following equality

$$b_j \cdot (c_{jkj} \cdot d_{jk}) = (b_j \cdot c_{jkj}) \cdot d_{jk} – (c_{jkj} \cdot b_j) \cdot d_{jk} + c_{jkj} \cdot (b_j \cdot d_{jk}).$$

Since $C$ is an ideal, $b_j \cdot c_{jkj}, c_{jkj} \cdot b_j \in C$, so that $(b_j \cdot c_{jkj}) \cdot d_{jk}, (c_{jkj} \cdot b_j) \cdot d_{jk} \in C \bullet D$. Since $D$ is an ideal, $b_j \cdot d_{jk} \in D$, so that $c_{jkj} \cdot (b_j \cdot d_{jk}) \in C \bullet D$, and hence $b_j \cdot (c_{jkj} \cdot d_{jk}) \in C \bullet D$. Thus $b_j \cdot (c_{jkj} \cdot d_{jk}) \in (B \bullet D) \cap (C \bullet D)$. Since this is true for all indexes j, k, we obtain that $x \in (B \bullet D) \cap (C \bullet D)$.

**2.5. LEMMA.** *Let $L$ be a pre – Lie algebra, B, C are ideals of L such that $L = B + C$. If $D = B \cap C$, then $L \bullet D = B \bullet D + C \bullet D$, $D \bullet L = D \bullet B + D \bullet C$.*

**PROOF.** Let $x$ be an arbitrary element of $L$, and let $d$ be the arbitrary elements of $D$, then $x = b + c$ where $b \in B$, $c \in C$. We have

$$x \cdot d = (b + c) \cdot d = b \cdot d + c \cdot d \in B \bullet D + C \bullet D,$$
$$d \cdot x = d \cdot (b + c) = d \cdot b + d \cdot c \in D \bullet B + D \bullet C,$$

so that $L \bullet D \leq B \bullet D + C \bullet D$ and $D \bullet L \leq D \bullet B + D \bullet C$. The converse inclusion is obvious.

**2.6. LEMMA.** *Let $L$ be a pre – Lie algebra, B, C are ideals of L, $D = B \cap C$. If $L = B + C$, then $(B \bullet D) \cap (C \bullet D)$ is a left ideal of L.*

**PROOF.** Let $x$ be an arbitrary element of $L$, then $x = b + c$ where $b \in B$, $c \in C$. Let $u$ be an arbitrary element of $(B \bullet D) \cap (C \bullet D)$. We have $u = \sum_{1 \leq j \leq n} b_j \cdot d_j$ for some elements $b_j \in B$, $d_j \in D$, $1 \leq j \leq n$. Then

$$x \cdot u = (b + c) \cdot u = (b + c) \cdot u = b \cdot u + c \cdot u =$$

$$b \cdot u + c \cdot \sum_{1 \le j \le n} b_j \cdot d_j = b \cdot u + \sum_{1 \le j \le n} c \cdot (b_j \cdot d_j).$$

Equality

$$(c \cdot b_j) \cdot d_j - c \cdot (b_j \cdot d_j) = (b_j \cdot c) \cdot d_j - b_j \cdot (c \cdot d_j)$$

implies that

$$c \cdot (b_j \cdot d_j) = (c \cdot b_j) \cdot d_j - (b_j \cdot c) \cdot d_j + b_j \cdot (c \cdot d_j).$$

Since B is an ideal of L, then $c \cdot b_j$, $b_j \cdot c \in B$, so that $(c \cdot b_j) \cdot d_j$, $(b_j \cdot c) \cdot d_j \in B \bullet D$. The fact that D is an ideal of L implies that $c \cdot d_j \in D$, so that $b_j \cdot (c \cdot d_j) \in B \bullet D$. Since $u \in D$, $b \cdot u \in B \bullet D$. It follows that $x \cdot u \in B \bullet D$. Using similar arguments, we obtain that $x \cdot u \in C \bullet D$, so that $x \cdot u \in (B \bullet D) \cap (C \bullet D)$.

**2.7. THEOREM.** *Let L be a pre – Lie algebra, B, C be ideals of L such that L = B + C. If $B^3 = C^3 = <0>$, then $L^6 = <0>$.*

**PROOF.** Put $D = B \cap C$. Using **Corollary 1.5** we obtain that $(L/D)^3 = <0>$. The equality $(L/D)^3 = (L^3 + D)/D$ implies that $L^3 \le D$. It follows that $L^4 = L \bullet L^3 \le L \bullet D$. **Lemma 2.5** implies that $L \bullet D = B \bullet D + C \bullet D$. This equality shows that every element of $L \bullet D$ is a sum of elements of one of the following forms $b \cdot d$ or $c \cdot d$ where $b \in B$, $c \in C$, $d \in D$.

Let x be an arbitrary element of L, then $x = b_1 + c_1$ where $b_1 \in B$, $c_1 \in C$. We have

$$x \cdot (b \cdot d) = (b_1 + c_1) \cdot (b \cdot d) = b_1 \cdot (b \cdot d) + c_1 \cdot (b \cdot d),$$
$$x \cdot (c \cdot d) = (b_1 + c_1) \cdot (c \cdot d) = b_1 \cdot (c \cdot d) + c_1 \cdot (c \cdot d).$$

The inclusion $D \le B$ implies that $b \cdot d \in B^2$, so that $b_1 \cdot (b \cdot d) \in B \bullet B^2 = B^3 = <0>$. Similarly, $c_1 \cdot (c \cdot d) = 0$. Thus $L \bullet (L \bullet D)$ as a subspace is generated by elements of one of the following forms $c_1 \cdot (b \cdot d)$ or $b_1 \cdot (c \cdot d)$.

Applying the above arguments, we obtain that $L \bullet (L \bullet (L \bullet D))$ as a subspace is generated the elements of one of the following forms

$$c_2 \cdot (c_1 \cdot (b \cdot d)), c_2 \cdot (b_1 \cdot (c \cdot d)), b_2 \cdot (c_1 \cdot (b \cdot d)), b_2 \cdot (b_1 \cdot (c \cdot d)).$$

Since D is an ideal, $b \cdot d, c \cdot d \in D = B \cap C$, so that $c_2 \cdot (c_1 \cdot (b \cdot d)) \in C \bullet C^2 = C^3 = <0>$, $b_2 \cdot (b_1 \cdot (c \cdot d)) \in B \bullet B^2 = B^3 = <0>$.

The equality $(c_1 \cdot b) \cdot d - c_1 \cdot (b \cdot d) = (b \cdot c_1) \cdot d - b \cdot (c_1 \cdot d)$ implies that

$$c_1 \cdot (b \cdot d) = (c_1 \cdot b) \cdot d - (b \cdot c_1) \cdot d + b \cdot (c_1 \cdot d).$$

Since B is an ideal, $c_1 \cdot b, b \cdot c_1 \in B$, so that $(c_1 \cdot b) \cdot d, (b \cdot c_1) \cdot d \in B \bullet D$. Since D is an ideal, $c_1 \cdot d \in D$, so that $b \cdot (c_1 \cdot d) \in B \bullet D$. Thus, $c_1 \cdot (b \cdot d) \in B \bullet D \le B \bullet B = B^2$. Then $b_2 \cdot (c_1 \cdot (b \cdot d)) \in B \bullet B^2 = B^3 = <0>$. Similarly, $b_1 \cdot (c \cdot d) \in C \bullet D \le C \bullet C = C^2$. Then $c_2 \cdot (b_1 \cdot (c \cdot d)) \in C \bullet C^2 = C^3 = <0>$. Hence $L \bullet (L \bullet (L \bullet D)) = <0>$.

The inclusion $L^4 \le L \bullet D$ implies that

$$L^5 = L \bullet L^4 \le L \bullet (L \bullet D), \; L^6 = L \bullet L^5 \le L \bullet (L \bullet (L \bullet D)) = <0>.$$

**2.8. LEMMA.** *Let L be a pre – Lie algebra, B, C be the ideals of L, $D = B \cap C$. Suppose that $C^3 = <0>$. Then $c \cdot (b_j \cdot (b_{j-1} \cdot (b_{j-2} \cdot \ldots \cdot (b_1 \cdot (d_1 \cdot (v_{m-1} \cdot \ldots (v_2 \cdot (v_1 \cdot d)) \ldots ) = 0$ for all positive integer j, m, j > m where $c \in C$, $b_j, b_{j-1}, \ldots b_1, v_m, v_{m-1}, \ldots v_1 \in B$, $d_1, d \in D$.*

**PROOF.** Let $u = c \cdot (b_j \cdot (b_{j-1} \cdot (b_{j-2} \cdot \ \ldots \ \cdot (b_1 \cdot \ (d_1 \cdot (v_{m-1} \cdot \ldots \ (v_2 \cdot \ (v_1 \cdot d)) \ldots)$. First, we will apply induction on $j$. If $j = 0$, then $u = c \cdot (d_1 \cdot (v_{m-1} \cdot \ldots \ (v_2 \cdot \ (v_1 \cdot d)) \ldots)$. Since $D$ is an ideal, then $(v_{m-1} \cdot \ldots \ (v_2 \cdot \ (v_1 \cdot d)) \ldots) = d_2 \in D$. The inclusion $D \leq C$ implies that $u = c \cdot (d_1 \cdot (v_{m-1} \cdot \ldots \ (v_2 \cdot \ (v_1 \cdot d)) \ldots) = c \cdot (d_1 \cdot d_2) \in C^3 = <0>$.

Let $j > 0$ and put $d_3 = (b_{j-1} \cdot (b_{j-2} \cdot \ \ldots \ \cdot (b_1 \cdot \ (d_1 \cdot (v_{m-1} \cdot \ldots \ (v_2 \cdot \ (v_1 \cdot d)) \ldots)$. Then $u = c \cdot (b_j \cdot (b_{j-1} \cdot (b_{j-2} \cdot \ \ldots \ \cdot (b_1 \cdot \ (d_1 \cdot (v_{m-1} \cdot \ldots \ (v_2 \cdot \ (v_1 \cdot d)) \ldots) = c \cdot (b_j \cdot d_3)$. We have $c \cdot (b_j \cdot d_3) = (c \cdot b_j) \cdot d_3 - (b_j \cdot c) \cdot d_3 + b_j \cdot (c \cdot d_3)$. Since $C$ is an ideal, $c \cdot b_j$, $b_j \cdot c \in C$, and then

$$(c \cdot b_j) \cdot d_3 = (c \cdot b_j) \cdot (b_{j-1} \cdot (b_{j-2} \cdot \ \ldots \ \cdot (b_1 \cdot \ (d_1 \cdot (v_{m-1} \cdot \ldots \ (v_2 \cdot \ (v_1 \cdot d)) \ldots).$$

Using the induction hypothesis we obtain that $(c \cdot b_j) \cdot d_3 = 0$. Similarly we can prove that $(b_j \cdot c) \cdot d_3 = 0$. Also we have

$$(c \cdot d_3) = (c \cdot (b_{j-1} \cdot (b_{j-2} \cdot \ \ldots \ \cdot (b_1 \cdot \ (d_1 \cdot (v_{m-1} \cdot \ldots \ (v_2 \cdot \ (v_1 \cdot d)) \ldots),$$

and using again the induction hypothesis we obtain that $(c \cdot d_3) = 0$, and hence $b_j \cdot (c \cdot d_3) = 0$. Thus

$$\mathbf{c \cdot (b_j \cdot (b_{j-1} \cdot (b_{j-2} \cdot \ \ldots \ \cdot (b_1 \cdot \ (d_1 \cdot (v_{m-1} \cdot \ldots \ (v_2 \cdot \ (v_1 \cdot d)) \ldots) = 0.}$$

Let $L$ be a pre – Lie algebra, $J, K$ be the subalgebras of $L$. Put $\lambda_1(J, K) = J \cdot K$, $\lambda_2(J, K) = J \cdot (J \cdot K) = J \cdot \lambda_1(J, K)$, and recursively $\lambda_{n+1}(J, K) = J \cdot \lambda_n(J, K)$ for all positive integers $n$.

**2.9. PROPOSITION.** *Let $L$ be a pre – Lie algebra, $B$ and $C$ be ideals of $L$, $D = B \cap C$. Suppose that $C^3 = <0>$. Then*

$$\lambda_k(L, D) \leq \lambda_k(B, D) + C \bullet \lambda_{k-1}(B, D) + \lambda_{k-2}(B, D \cdot D) + \lambda_{k-3}(B, D \cdot \lambda_1(B, D)) + \lambda_{k-4}(B, D \cdot \lambda_2(B, D)) + \ldots + \lambda_2(B, D \cdot \lambda_{k-4}(B, D)) + \lambda_1(B, D \cdot \lambda_{k-3}(B, D)) + D \cdot \lambda_{k-2}(B, D))$$

*for all positive integers $k$.*

**PROOF. Lemma 1.11** implies that $L \bullet D = B \bullet D + C \bullet D$. This equality shows that every element of $L \bullet D$ is a sum of the elements of one of the following forms $b \cdot d$ or $c \cdot d$ where $b \in B, c \in C, d \in D$.

Let $x$ be the arbitrary element of $L$. Then $x = b_1 + c_1$ where $b_1 \in B, c_1 \in C$. We have

$$x \cdot (b \cdot d) = (b_1 + c_1) \cdot (b \cdot d) = b_1 \cdot (b \cdot d) + c_1 \cdot (b \cdot d),$$
$$x \cdot (c \cdot d) = (b_1 + c_1) \cdot (c \cdot d) = b_1 \cdot (c \cdot d) + c_1 \cdot (c \cdot d).$$

The inclusion $D \leq C$ implies that $c \cdot d \in C^2$, so that $c_1 \cdot (c \cdot d) \in C \bullet C^2 = C^3 = <0>$.Thus, $L \bullet (L \bullet D)$ as a subspace is generated by elements of one of the following forms $b_1 \cdot (b \cdot d)$, $c_1 \cdot (b \cdot d)$, $b_1 \cdot (c \cdot d)$. We have

$$b_1 \cdot (c \cdot d) = (b_1 \cdot c) \cdot d - (c \cdot b_1) \cdot d + c \cdot (b_1 \cdot d).$$

Note that $b_1 \cdot c, c \cdot b_1 \in B \cap C = D$, so that $b_1 \cdot (c \cdot d) \in D \bullet D + C \bullet (B \bullet D) \leq D^2 + C \bullet B^2 \leq B^2 + C \bullet B^2$. Thus, we see that

$$\lambda_2(L, D) \leq \lambda_2(B, D) + C \bullet \lambda_1(B, D) + D \bullet D.$$

In other words,

$L \bullet (L \bullet D)$ is a subspace of the linear envelope of a set consisting of elements of one of the following forms:

$$d_1 \cdot d,\ b_1 \cdot (b \cdot d),\ c \cdot (b \cdot d),\ \ b, b_1 \in B, c \in C, d \in D.$$

It is not difficult to show that $\lambda_3(L, D) = L \bullet (L \bullet (L \bullet D))$ is a subspace of the linear envelop of a set consisting of elements of one of the following forms:

$$b_2 \cdot (d_1 \cdot d),\ b_2 \cdot (b_1 \cdot (b \cdot d)),\ b_2 \cdot (c \cdot (b \cdot d)),\ c \cdot (d_1 \cdot d),\ c \cdot (b_1 \cdot (b \cdot d)),\ c_1 \cdot (c \cdot (b \cdot d)),$$
$$b, b_1, b_2 \in B,\ c, c_1 \in C,\ d, d_1 \in D.$$

The inclusion $D \le C$ and the fact that $D$ is an ideal imply that $c_1 \cdot (c \cdot (b \cdot d)) \in C^3 = \langle 0 \rangle$. Similarly, $c \cdot (d_1 \cdot d) \in C^3 = \langle 0 \rangle$.

Consider element $b_2 \cdot (c \cdot (b \cdot d))$. Put $d_1 = (b \cdot d))$. Since $D$ is an ideal, $d_1 \in D$. We have

$$b_2 \cdot (c \cdot (b \cdot d)) = b_2 \cdot (c \cdot d_1) = (b_2 \cdot c) \cdot d_1 - (c \cdot b_2) \cdot d_1 + c \cdot (b_2 \cdot d_1) =$$
$$(b_2 \cdot c) \cdot (b \cdot d) - (c \cdot b_2) \cdot (b \cdot d)) + c \cdot (b_2 \cdot (b \cdot d)).$$

Since $b_1 \cdot c,\ c \cdot b_1 \in B \cap C = D$, $b_2 \cdot (c \cdot (b \cdot d)) \in D \bullet (B \bullet D) + C \bullet \lambda_2(B, D) \le D \bullet \lambda_1(B, D) + C \bullet \lambda_2(B, D)$. Thus we see that

$$\lambda_3(L, D) \le \lambda_3(B, D) + C \bullet \lambda_2(B, D) + B \bullet (D \bullet D) + D \bullet (B \bullet D).$$

In other words, $\lambda_3(L, D)$ is a subspace of the linear envelope of a set consisting of elements of one of the following forms:

$$b_2 \cdot (d_1 \cdot d),\ b_2 \cdot (b_1 \cdot (b \cdot d)),\ d_1 \cdot (b \cdot d),\ c \cdot (b_1 \cdot (b \cdot d)),\ b, b_1, b_2 \in B,\ c, c_1 \in C,\ d, d_1 \in D.$$

It is not difficult to show that $\lambda_4(B, D) = L \bullet (L \bullet (L \bullet (L \bullet D)))$ is a subspace of the linear envelope of a set consisting of elements of one of the following forms:

$$b_1 \cdot (b \cdot (d_1 \cdot d)),\ b_3 \cdot (b_2 \cdot (b_1 \cdot (b \cdot d))),\ b_1 \cdot (d_1 \cdot (b \cdot d)),\ b_2 \cdot (c \cdot (b_1 \cdot (b \cdot d))),$$
$$c \cdot (b_2 \cdot (d_1 \cdot d)),\ c \cdot (b_2 \cdot (b_1 \cdot (b \cdot d)),\ c \cdot (d_1 \cdot (b \cdot d)),\ c_1 \cdot (c \cdot (b_1 \cdot (b \cdot d)),$$
$$b, b_1, b_2, b_3 \in B,\ c, c_1 \in C,\ d, d_1 \in D.$$

Consider the element $b_2 \cdot (c \cdot (b_1 \cdot (b \cdot d)))$. Put $d_1 = b_1 \cdot (b \cdot d)$. Since $D$ is an ideal, $d_1 \in D$. We have

$$b_2 \cdot (c \cdot (b_1 \cdot (b \cdot d))) = b_2 \cdot (c \cdot d_1) = (b_2 \cdot c) \cdot d_1 - (c \cdot b_2) \cdot d_1 + c \cdot (b_2 \cdot d_1),$$
$$(b_2 \cdot c) \cdot d_1 = (b_2 \cdot c) \cdot (b_1 \cdot (b \cdot d)) \in D \bullet \lambda_2(B, D),\ (c \cdot b_2) \cdot d_1 \in D \bullet \lambda_2(B, D),$$
$$c \cdot (b_2 \cdot d_1) = c \cdot (b_2 \cdot (b_1 \cdot (b \cdot d))) \in C \bullet \lambda_3(B, D).$$

Consider the element $c \cdot (b_2 \cdot (d_1 \cdot d))$. Put $d_2 = d_1 \cdot d$. We have

$$c \cdot (b_2 \cdot (d_1 \cdot d)) = c \cdot (b_2 \cdot d_2) = (c \cdot b_2) \cdot d_2 - (b_2 \cdot c) \cdot d_2 + b_2 \cdot (c \cdot d_2),$$
$$(c \cdot b_2) \cdot d_2 = (c \cdot b_2) \cdot (d_1 \cdot d) \in D^3 \le C^3 = \langle 0 \rangle,\ (b_2 \cdot c) \cdot d_2 = 0,\ c \cdot d_2 = c \cdot (d_1 \cdot d) \in C^3 = \langle 0 \rangle,$$

so that $c \cdot (b_2 \cdot (d_1 \cdot d)) = 0$. We have $c \cdot (d_1 \cdot (b \cdot d)) \in C \bullet D^2 \le C^3 = \langle 0 \rangle$, so that $c \cdot (d_1 \cdot (b \cdot d)) = 0$.

Since $b_1 \cdot (b \cdot d)) \in D$, $c_1 \cdot (c \cdot (b_1 \cdot (b \cdot d)) \in C \bullet (C \bullet D) \le C^3 = \langle 0 \rangle$, so that $c_1 \cdot (c \cdot (b_1 \cdot (b \cdot d)) = 0$.

Thus, we see that

$\lambda_4(L, D) \leq \lambda_4(B, D) + C \cdot \lambda_3(B, D) + B \cdot (B \cdot (D \cdot D)) + B \cdot (D \cdot (B \cdot D)) + D \cdot (B \cdot (B \cdot D) =$
$\lambda_4(B, D) + C \cdot \lambda_3(B, D) + D \cdot \lambda_2(B, D) + B \cdot (D \cdot \lambda_1(B, D)) + \lambda_2(B, D \cdot D)$

Let k be a positive integer, and suppose that we have already proved that

$\lambda_j(L, D) \leq \lambda_j(B, D) + C \cdot \lambda_{j-1}(B, D) + \lambda_{j-2}(B, D \cdot D) + \lambda_{j-3}(B, D \cdot \lambda_1(B, D)) + \lambda_{j-4}(B, D \cdot \lambda_2(B, D)) + \ldots + \lambda_2(B, D \cdot \lambda_{j-4}(B, D)) + \lambda_1(B, D \cdot \lambda_{j-3}(B, D)) + D \cdot \lambda_{j-2}(B, D))$

for all $j < k$.

We have $\lambda_k(L, D) = L \cdot \lambda_{k-1}(L, D)$

$\lambda_k(L, D) = L \cdot \lambda_{k-1}(L, D) = (B + C) \cdot \lambda_{k-1}(L, D) = B \cdot \lambda_{k-1}(L, D) + C \cdot \lambda_{k-1}(L, D),$

$B \cdot \lambda_{k-1}(L, D) \leq B \cdot (\lambda_{k-1}(B, D) + C \cdot \lambda_{k-2}(B, D) + \lambda_{k-3}(B, D \cdot D) + \lambda_{k-4}(B, D \cdot \lambda_1(B, D)) +$
$+ \ldots + \lambda_2(B, D \cdot \lambda_{k-5}(B, D)) + \lambda_1(B, D \cdot \lambda_{j-4}(B, D)) + D \cdot \lambda_{j-3}(B, D))) =$
$B \cdot (\lambda_{k-1}(B, D)) + B \cdot (C \cdot \lambda_{k-2}(B, D)) + B \cdot (\lambda_{k-3}(B, D \cdot D)) + B \cdot (\lambda_{k-4}(B, D \cdot \lambda_1(B, D))) +$
$+ \ldots + B \cdot (\lambda_2(B, D \cdot \lambda_{k-5}(B, D))) + B \cdot (\lambda_1(B, D \cdot \lambda_{k-4}(B, D))) + B \cdot (D \cdot \lambda_{k-3}(B, D)))) =$
$\lambda_k(B, D)) + B \cdot (C \cdot \lambda_{k-2}(B, D)) + \lambda_{k-2}(B, D \cdot D)) + \lambda_{k-3}(B, D \cdot \lambda_1(B, D))) + \ldots +$
$\lambda_3(B, D \cdot \lambda_{k-5}(B, D))) + \lambda_2(B, D \cdot \lambda_{k-4}(B, D))) + \lambda_1(B, D \cdot \lambda_{k-3}(B, D)).$

The subspace $B \cdot (C \cdot \lambda_{k-2}(B, D))$ is a linear envelope of the set of all elements having the form $b \cdot (c \cdot x)$ where $b \in B, c \in B, x \in \lambda_{k-2}(B, D)$. We have

$$b \cdot (c \cdot x) = (b \cdot c) \cdot x - (c \cdot b) \cdot x + c \cdot (b \cdot x)$$

Since $B, C$ are ideals, $c \cdot b, b \cdot c \in B \cap C = D$, so that $c \cdot b, b \cdot c \in D$ and $(c \cdot b) \cdot x$, $(b \cdot c) \cdot x \in D \cdot \lambda_{k-2}(B, D))$. Since $x \in \lambda_{k-2}(B, D)$, $b \cdot x \in \lambda_{k-1}(B, D))$ and $c \cdot (b \cdot x) \in C \cdot \lambda_{k-1}(B, D))$. Thus $B \cdot (C \cdot \lambda_{k-2}(B, D)) \leq D \cdot \lambda_{k-2}(B, D)) + C \cdot \lambda_{k-1}(B, D))$.

Furthermore,

$C \cdot \lambda_{k-1}(L, D) \leq \mathbf{C} \cdot \mathbf{(}\lambda_{k-1}(B, D) + C \cdot \lambda_{k-2}(B, D) + \lambda_{k-3}(B, D \cdot D) + \lambda_{k-4}(B, D \cdot \lambda_1(B, D)) +$
$+ \ldots + \lambda_2(B, D \cdot \lambda_{k-5}(B, D)) + \lambda_1(B, D \cdot \lambda_{j-4}(B, D)) + D \cdot \lambda_{j-3}(B, D))) =$
$\mathbf{C} \cdot \mathbf{(}\lambda_{k-1}(B, D)) + \mathbf{C} \cdot (C \cdot \lambda_{k-2}(B, D)) + \mathbf{C} \cdot \mathbf{(}\lambda_{k-3}(B, D \cdot D)) + \mathbf{C} \cdot \mathbf{(}\lambda_{k-4}(B, D \cdot \lambda_1(B, D))) +$
$+ \ldots + \mathbf{C} \cdot \mathbf{(}\lambda_2(B, D \cdot \lambda_{k-5}(B, D))) + \mathbf{C} \cdot \mathbf{(}\lambda_1(B, D \cdot \lambda_{k-4}(B, D))) + \mathbf{C} \cdot \mathbf{(}D \cdot \lambda_{k-3}(B, D)))).$

Since $D$ is an ideal, $\lambda_{k-2}(B, D)) \leq D \leq C$, so that $\mathbf{C} \cdot \mathbf{(}C \cdot \lambda_{k-2}(B, D)) \leq C^3 = \langle 0 \rangle$. Using **Lemma 2.8** we obtain that every term of this sum is zero. Hence

$\lambda_k(L, D) \leq \lambda_k(B, D) + C \cdot \lambda_{k-1}(B, D) + \lambda_{k-2}(B, D \cdot D) + \lambda_{k-3}(B, D \cdot \lambda_1(B, D)) + \lambda_{k-4}(B, D \cdot \lambda_2(B, D)) + \ldots + \lambda_2(B, D \cdot \lambda_{k-4}(B, D)) + \lambda_1(B, D \cdot \lambda_{k-3}(B, D)) +$
$D \cdot \lambda_{k-2}(B, D)).$

**THEOREM A.** *Let $L$ be a pre – Lie algebra, $B$, $C$ be ideals of $L$ such that $L = B + C$. Suppose that $B^k = C^3 = \langle 0 \rangle$. If $B$ is abelian, then $L^4 = \langle 0 \rangle$. If $k > 2$, then $L^{2k} = \langle 0 \rangle$.*

**PROOF.** Put $D = B \cap C$. If $B$ is abelian, then **Lemma 2.5** implies that $L \cdot D = B \cdot D + C \cdot D = C \cdot D$. It follows that

$$L \cdot (L \cdot D) = (B + C) \cdot (C \cdot D) = B \cdot (C \cdot D) + C \cdot (C \cdot D).$$

Since $D$ is an ideal, $C \bullet D \leq D \leq B$, so that $B \bullet (C \bullet D) = <0>$. The inclusion $D \leq C$ implies that $C \bullet (C \bullet D) \leq C \bullet (C \bullet C) = <0>$. Thus $L \bullet (L \bullet D) = <0>$. Using **Corollary 1.5** we obtain that $(L/D)^3 = <0>$. The equality $(L/D)^3 = (L^3 + D)/D$ implies that $L^3 \leq D$. It follows that

$$L^4 = L \bullet L^3 \leq L \bullet D, L^5 = L \bullet L^4 \leq L \bullet (L \bullet D) = <0>$$

If $B^3 = <0>$, then **Theorem 2.7** implies that $L^6 = <0>$.

Suppose now that $k > 3$. Using **Corollary 1.5** we obtain that $(L/D)^k = <0>$. The equality $(L/D)^k = (L^k + D)/D$ implies that $L^k \leq D$. By **Proposition 2.9**

$$\lambda_k(L, D) \leq \lambda_k(B, D) + C \bullet \lambda_{k-1}(B, D) + \lambda_{k-2}(B, D \cdot D) + \lambda_{k-3}(B, D \cdot \lambda_1(B, D)) + \lambda_{k-4}(B, D \cdot \lambda_2(B, D)) + \ldots + \lambda_2(B, D \cdot \lambda_{k-4}(B, D)) + \lambda_1(B, D \cdot \lambda_{k-3}(B, D)) + D \cdot \lambda_{k-2}(B, D)).$$

The inclusion $D \leq B$ implies that $B^k$ includes $\lambda_{k-2}(B, D \cdot D)$, $\lambda_{k-3}(B, D \cdot \lambda_1(B, D))$, $\lambda_{k-4}(B, D \cdot \lambda_2(B, D))$, **. . .** , $\lambda_2(B, D \cdot \lambda_{k-4}(B, D))$, $\lambda_1(B, D \cdot \lambda_{k-3}(B, D))$, $D \cdot \lambda_{k-2}(B, D))$.
Also we have $\lambda_k(B, D) \leq B^{k+1}$, $\lambda_{k-1}(B, D) \leq B^k$ and $C \bullet \lambda_{k-1}(B, D) \leq C \bullet B^k$. Therefore, the equality $B^k = <0>$ implies that $\lambda_k(L, D) = <0>$. The inclusion $L^k \leq D$ implies that $L^{k+k} \leq \lambda_k(L, D)$, and we obtain that $L^{2k} = <0>$.

**2.10. LEMMA.** *Let $L$ be a pre – Lie algebra, $B$, $C$ be ideals of $L$ such that $L = B + C$, $D = B \cap C$. If $B^{(3)} = C^{(3)} = <0>$, then $J = (D \bullet C) \bullet B + (D \bullet B) \bullet C + (B \bullet D) \bullet C + (C \bullet D) \bullet B$ is an abelian ideal of $L$. Moreover, $J \bullet D = <0>$.*

**PROOF.** It is not hard to see that $J$ is a linear envelope of the subset

$$\mathbf{M} = \{(d \cdot b) \cdot c, (d \cdot c) \cdot b, (b \cdot d) \cdot c \text{ or } (c \cdot d) \cdot b \mid b \in B, c \in C, d \in D \}.$$

Let $x$ be an arbitrary element of $L$. Then $x = b_1 + c_1$ where $b_1 \in B$, $c_1 \in C$. If $y$ is an arbitrary element of $M$, then

$$y \cdot x = y \cdot (b_1 + c_1) = y \cdot b_1 + y \cdot c_1, x \cdot y = (b_1 + c_1) \cdot y = b_1 \cdot y + c_1 \cdot y.$$

Thus, it is enough to prove that $y \cdot b_1$, $y \cdot c_1$, $b_1 \cdot y$, $c_1 \cdot y \in J$ for every element $y \in M$. We have

$$((d \cdot b) \cdot c) \cdot c_1, ((b \cdot d) \cdot c) \cdot c_1 \in (D \bullet C) \bullet C \leq (C \bullet C) \bullet C = C^{(3)} = <0>,$$
$$((d \cdot c) \cdot b) \cdot b_1, ((c \cdot d) \cdot b) \cdot b_1 \in (D \bullet B) \bullet B \leq (B \bullet B) \bullet B = B^{(3)} = <0>,$$
$$((d \cdot b) \cdot c) \cdot b_1 \in (D \bullet C) \bullet B \leq J, ((d \cdot c) \cdot b) \cdot c_1 \in (D \bullet B) \bullet C \leq J,$$
$$((b \cdot d) \cdot c) \cdot b_1 \in (D \bullet C) \bullet B \leq J, ((c \cdot d) \cdot b) \cdot c_1 \in (D \bullet B) \bullet C \leq J.$$

On the other hand,

$$b_1 \cdot ((d \cdot b) \cdot c) = (d \cdot b) \cdot (b_1 \cdot c) – ((d \cdot b) \cdot b_1) \cdot c + (b_1 \cdot (d \cdot b)) \cdot c,$$

$(d \cdot b) \cdot b_1 \in (D \bullet B) \bullet B \leq (B \bullet B) \bullet B = B^{(3)} = <0>$, so that $((d \cdot b) \cdot b_1) \cdot c = 0$,
Since $C$, $D$ are ideals, $b_1 \cdot c \in C$, $d \cdot b \in D$, so that $(d \cdot b) \cdot (b_1 \cdot c) \in (D \bullet B) \bullet C$, $(b_1 \cdot (d \cdot b)) \cdot c \in (B \bullet D) \bullet C$, and hence $b_1 \cdot ((d \cdot b) \cdot c) \in (D \bullet B) \bullet C + (B \bullet D) \bullet C \leq J$;

$$c_1 \cdot ((d \cdot b) \cdot c) = (c_1 \cdot (d \cdot b)) \cdot c + (d \cdot b) \cdot (c \cdot c_1) – ((d \cdot b) \cdot c_1) \cdot c,$$

$((d \cdot b) \cdot c_1) \cdot c \in (D \bullet C) \bullet C \leq (C \bullet C) \bullet C = C^{(3)} = <0>$,
$(c_1 \cdot (d \cdot b)) \cdot c \in (C \bullet D) \bullet C \leq (C \bullet C) \bullet C = C^{(3)} = <0>$,
$(d \cdot b) \cdot (c \cdot c_1) \in (D \bullet B) \bullet C$, and hence $c_1 \cdot ((d \cdot b) \cdot c) \in (D \bullet B) \bullet C \leq J$;

$$b_1 \cdot ((b \cdot d) \cdot c) = (b_1 \cdot (b \cdot d)) \cdot c + (b \cdot d) \cdot (b_1 \cdot c) - ((b \cdot d) \cdot b_1) \cdot c,$$

$(b \cdot d) \cdot b_1 \in (B \bullet D) \bullet B \le (B \bullet B) \bullet B = B^{(3)} = \langle 0\rangle$, so that $((b \cdot d) \cdot b_1) \cdot c = 0$,
$(b \cdot d) \cdot (b_1 \cdot c) \in (B \bullet D) \bullet C$, $(b_1 \cdot (b \cdot d)) \cdot c \in (B \bullet D) \bullet C$
and hence $\mathbf{b_1 \cdot ((b \cdot d) \cdot c)} \in (B \bullet D) \bullet C \le J$;

$$c_1 \cdot ((b \cdot d) \cdot c) = (c_1 \cdot (b \cdot d)) \cdot c + (b \cdot d) \cdot (c_1 \cdot c) - ((b \cdot d) \cdot c_1) \cdot c,$$

$(c_1 \cdot (b \cdot d)) \cdot c \in (C \bullet D) \bullet C \le (C \bullet C) \bullet C = C^{(3)} = \langle 0\rangle$,
$((b \cdot d) \cdot c_1) \cdot c \in (D \bullet C) \bullet C \le (C \bullet C) \bullet C = C^{(3)} = \langle 0\rangle$,
$(b \cdot d) \cdot (c_1 \cdot c) \in (B \bullet D) \bullet C$ and hence $c_1 \cdot ((b \cdot d) \cdot c) = (b \cdot d) \cdot (c_1 \cdot c) \in (B \bullet D) \bullet C \le J$;

$$b_1 \cdot ((d \cdot c) \cdot b) = (b_1 \cdot (d \cdot c)) \cdot b + (d \cdot c) \cdot (b_1 \cdot b) - ((d \cdot c) \cdot b_1) \cdot b,$$

$(b_1 \cdot (d \cdot c)) \cdot b \in (B \bullet D) \bullet B \le (B \bullet B) \bullet B = B^{(3)} = \langle 0\rangle$,
$((d \cdot c) \cdot b_1) \cdot b \in (D \bullet B) \bullet B \le (B \bullet B) \bullet B = B^{(3)} = \langle 0\rangle$,
$(d \cdot c) \cdot (b_1 \cdot b) \in (D \bullet C) \bullet B$ and hence $b_1 \cdot ((d \cdot c) \cdot b) = (d \cdot c) \cdot (b_1 \cdot b) \in (D \bullet C) \bullet B \le J$;

$$c_1 \cdot ((d \cdot c) \cdot b) = (c_1 \cdot (d \cdot c)) \cdot b + (d \cdot c) \cdot (c_1 \cdot b) - ((d \cdot c) \cdot c_1) \cdot b,$$

$(d \cdot c) \cdot c_1 \in (D \bullet C) \bullet C \le (C \bullet C) \bullet C = C^{(3)} = \langle 0\rangle$, so that $((d \cdot c) \cdot c_1) \cdot b = 0$,
$(d \cdot c) \cdot (c_1 \cdot b) \in (D \bullet C) \bullet B$, $(c_1 \cdot (d \cdot c)) \cdot b \in (C \bullet D) \bullet B$,
and hence $\mathbf{c_1 \cdot ((d \cdot c) \cdot b)} \in (C \bullet D) \bullet B + (D \bullet C) \bullet B \le J$;

$$b_1 \cdot ((c \cdot d) \cdot b) = (b_1 \cdot (c \cdot d)) \cdot b + (c \cdot d) \cdot (b_1 \cdot b) - ((c \cdot d) \cdot b_1) \cdot b,$$

$(b_1 \cdot (c \cdot d)) \cdot b \in (B \bullet D) \bullet B \le (B \bullet B) \bullet B = B^{(3)} = \langle 0\rangle$,
$((c \cdot d) \cdot b_1) \cdot b \in (D \bullet B) \bullet B \le (B \bullet B) \bullet B = B^{(3)} = \langle 0\rangle$,
$(c \cdot d) \cdot (b_1 \cdot b) \in (C \bullet D) \bullet B$, so that $b_1 \cdot ((c \cdot d) \cdot b) = (c \cdot d) \cdot (b_1 \cdot b) \in (C \bullet D) \bullet B \le J$;

$$c_1 \cdot ((c \cdot d) \cdot b) = (c_1 \cdot (c \cdot d)) \cdot b + (c \cdot d) \cdot (c_1 \cdot b) - ((c \cdot d) \cdot c_1) \cdot b,$$

$(c \cdot d) \cdot c_1 \in (C \bullet D) \bullet C \le (C \bullet C) \bullet C = C^{(3)} = \langle 0\rangle$, so that $((c \cdot d) \cdot c_1) \cdot b = 0$,
$(c \cdot d) \cdot (c_1 \cdot b) \in (C \bullet D) \bullet B$, $(c_1 \cdot (c \cdot d)) \cdot b \in (C \bullet D) \bullet B$,
and hence $c_1 \cdot ((c \cdot d) \cdot b) \in (C \bullet D) \bullet B \le J$.

Thus, we obtain that $J$ is an ideal of $L$. Furthermore, if $z$ is an arbitrary element of $D$, then

$$((d \cdot b) \cdot c) \cdot z \in (D \bullet C) \bullet D \le (C \bullet C) \bullet C = C^{(3)} = \langle 0\rangle,$$
$$((d \cdot c) \cdot b) \cdot z \in (D \bullet B) \bullet D \le (B \bullet B) \bullet B = B^{(3)} = \langle 0\rangle,$$
$$((c \cdot d) \cdot b) \cdot z \in (D \bullet B) \bullet D \le (B \bullet B) \bullet B = B^{(3)} = \langle 0\rangle,$$
$$((b \cdot d) \cdot c) \cdot z \in (D \bullet C) \bullet D \le (C \bullet C) \bullet C = C^{(3)} = \langle 0\rangle.$$

It follows that $J \bullet D = \langle 0\rangle$, in particular, the inclusion $J \le D$ implies that $J \bullet J = \langle 0\rangle$, so that $J$ is an abelian ideal.

**THEOREM B.** *Let $L$ be a pre – Lie algebra, $B$, $C$ be ideals of $L$ such that $L = B + C$. If $B^{(3)} = C^{(3)} = \langle 0\rangle$, then $B \cap C$ includes an abelian ideal $J$ such that $(L/J)^{(6)} = \langle 0\rangle$.*

**PROOF.** Put $D = B \cap C$. Using **Corollary 1.5** we obtain that $(L/D)^{(3)} = \langle 0\rangle$. The equality $(L/D)^{(3)} = (L^{(3)} + D)/D$ implies that $L^{(3)} \le D$. It follows that $L^{(4)} = L^{(3)} \bullet L \le$ $D \cdot L$. **Lemma 2.5** implies that $D \cdot L = D \cdot B + D \cdot C$. This equality shows that every element of $D \cdot L$ is a sum of elements of one of the following forms $d \cdot b$ or $d \cdot c$ where $b \in B, c \in C, d \in D$.

Let $x$ be the arbitrary element of $L$. Then $x = b_1 + c_1$ where $b_1 \in B, c_1 \in C$. We have

$$(d \cdot b) \cdot x = (d \cdot b) \cdot (b_1 + c_1) = (d \cdot b) \cdot b_1 + (d \cdot b) \cdot c_1,$$
$$(d \cdot c) \cdot x = (d \cdot c) \cdot (b_1 + c_1) = (d \cdot c) \cdot b_1 + (d \cdot c) \cdot c_1.$$

The inclusion $D \leq B$ implies that $d \cdot b \in B^{(2)}$, so that $(d \cdot b) \cdot b_1 \in B^{(2)} \cdot B = B^{(3)} = <0>$. Similarly, we obtain that $(d \cdot c) \cdot c_1 = 0$. Thus, $(D \cdot L) \cdot L$ as a subspace is generated by elements of one of the following forms $(d \cdot b) \cdot c_1$ or $(d \cdot c) \cdot b_1$.

Using the above arguments, we obtain that $((D \cdot L) \cdot L) \cdot L$ as a subspace is generated by elements of one of the following forms

$$((d \cdot b) \cdot c) \cdot b_1, ((d \cdot b) \cdot c) \cdot c_1, ((d \cdot c) \cdot b) \cdot b_1, ((d \cdot c) \cdot b) \cdot c_1,$$

where $b \in B, c \in C, d \in D$. Put $J = (D \cdot C) \cdot B + (D \cdot B) \cdot C + (B \cdot D) \cdot C + (C \cdot D) \cdot B$, then the above remark shows that $((D \cdot L) \cdot L) \cdot L \leq J$. **Lemma 2.10** implies that $J$ is an abelian ideal of $L$.

We have now $L^{(4)} = L^{(3)} \cdot L \leq D \cdot L$, $L^{(5)} = L^{(4)} \cdot L \leq (D \cdot L) \cdot L$, $L^{(6)} = L^{(5)} \cdot L \leq ((D \cdot L) \cdot L) \cdot L$. It follows that $L^{(6)} \leq J$, so that $(L/J)^{(6)} = <0>$.

**THEOREM C.** *Let $L$ be a pre – Lie algebra, $B, C$ be ideals of $L$ such that $L = B + C$, $B = \zeta_2(B)$, $C = \zeta_2(C)$. Then the following assertions hold*

*(i) $\zeta_1(B) \cap \zeta_1(C) \leq \zeta_1(L)$.*

*(ii) $B \cdot (C \cdot D) \leq \zeta_1(B) \cap \zeta_1(C)$, $C \cdot (B \cdot D) \leq \zeta_1(B) \cap \zeta_1(C)$.*

*(iii) $L \cdot D = B \cdot D + C \cdot D \leq (D \cap \zeta_1(B)) + (D \cap \zeta_1(C)) \leq \zeta_1(D)$,*

*(iv) $L \cdot (B \cdot D + C \cdot D) \leq B \cdot (C \cdot D) + C \cdot (B \cdot D) \leq \zeta_1(B) \cap \zeta_1(C) \leq \zeta_1(L)$, in particular, $B \cdot (C \cdot D) + C \cdot (B \cdot D)$ is an ideal of $L$.*

*(v) $L \cdot (L \cdot (L \cdot D) = <0>$.*

*(vi) $L^6 = <0>$.*

*(vii) $D \cdot L \leq D \cdot B + D \cdot C \leq (D \cap \zeta_1(B)) + (D \cap \zeta_1(C)) \leq \zeta_1(D)$,*

*(viii) $(D \cap \zeta_1(B)) + (D \cap \zeta_1(C))$ is an abelian ideal of $L$.*

*(ix) $\zeta_1(D)$ is an ideal of $L$.*

*(x) $(D \cdot L) \cdot L \leq (D \cdot C) \cdot B + (D \cdot B) \cdot C \leq (D \cap \zeta_1(B)) + (D \cap \zeta_1(C)) \leq \zeta_1(D)$, $((D \cdot L) \cdot L) \cdot L \leq (D \cdot C) \cdot B + (D \cdot B) \cdot C \leq \zeta_1(D)$.*

*(xi) $L^{(5)} \leq (D \cap \zeta_1(B)) + (D \cap \zeta_1(C)) \leq \zeta_1(D)$.*

*(xii) The factor – algebra $L/((D \cap \zeta_1(B)) + (D \cap \zeta_1(C)))$ is nilpotent and* ***uzl***$(A/((D \cap \zeta_1(\star, B)) + (D \cap \zeta_1(\star, C)))) \leq 3$.

**PROOF.** Put $D = B \cap C$, $B_1 = \zeta_1(B)$, $C_1 = \zeta_1(C)$.

(i) Let $x$ be an arbitrary element of $L$ and let $z$ be an arbitrary element of $\zeta_1(B) \cap \zeta_1(C)$. Then $x = b + c$ where $b \in B, c \in C$. We have

$$x \cdot z = (b + c) \cdot z = b \cdot z + c \cdot z, z \cdot x = z \cdot (b + c) = z \cdot b + z \cdot c.$$

Since $z \in \zeta_1(B)$, $b \cdot z = z \cdot b = 0$. Since $z \in \zeta_1(C)$, $c \cdot z = z \cdot c = 0$, and we obtain that $x \cdot z = z \cdot x = 0$. Hence $z \in \zeta_1(L)$.

(ii) Let $b$ be an arbitrary element of $B$, $c$ be an arbitrary element of $C$ and let $z$ be an arbitrary elements of $B \cap C$. Since $B$ is an ideal, $c \cdot z \in B$ and then $b \cdot (c \cdot z) \in B_1$.

On the other hand, $b \cdot (c \cdot z) = (b \cdot c) \cdot z + c \cdot (b \cdot z) - (c \cdot b) \cdot z$. Since $C$ is an ideal, $b \cdot c \in C$ and then $(b \cdot c) \cdot z \in C_1$. Since $D$ is an ideal, $b \cdot z \in D \leq C$ and then $c \cdot (b \cdot z) \in C_1$. Since $C$ is an ideal, $c \cdot b \in C$ and then $(c \cdot b) \cdot z \in C_1$. Hence $b \cdot (c \cdot z) \in C_1$. Similarly, $C \bullet (B \bullet D) \leq B_1 \cap C_1$.

(iii) The equality $L \bullet D = B \bullet D + C \bullet D$ was proven in **Lemma 2.5.**
Let $x = b + c$, $b \in B$, $c \in C$, be an arbitrary element of $A$ and let $d$ be an arbitrary element of $D$. As above, we have $x \cdot d = b \cdot d + c \cdot d$. Since $d \in B$, $b \cdot d \in B_1$. Since $d \in C$, $c \cdot d \in C_1$. On the other hand, $D$ is an ideal of $A$, and therefore, $b \cdot d \in B_1 \cap D$, $c \cdot d \in C_1 \cap D$, so that $L \bullet D \leq (D \cap B_1) + (D \cap C_1)$.

(iv) Let $x = b + c$, $b \in B$, $c \in C$, be an arbitrary element of $A$, $b_1$ be an arbitrary element of $B$, $c_1$ be anarbitrary element of $C$, and $d_1, d_2$ be arbitrary elements of $D$. We have

$$x \cdot (b_1 \cdot d_1) = (b + c) \cdot (b_1 \cdot d_1) = b \cdot (b_1 \cdot d_1) + c \cdot (b_1 \cdot d_1),$$
$$x \cdot (c_1 \cdot d_2) = (b + c) \cdot (c_1 \cdot d_2) = b \cdot (c_1 \cdot d_2) + c \cdot (c_1 \cdot d_2),$$
$$x \cdot (b_1 \cdot d_1 + c_1 \cdot d_2) = x \cdot (b_1 \cdot d_1) + x \cdot (c_1 \cdot d_2) =$$
$$b \cdot (b_1 \cdot d_1) + c \cdot (b_1 \cdot d_1) + b \cdot (c_1 \cdot d_2) + c \cdot (c_1 \cdot d_2).$$

Since $b_1, d \in B$, $b_1 \cdot d \in \zeta_1(B)$. It follows that $b \cdot (b_1 \cdot d) = 0$. Similarly, $c \cdot (c_1 \cdot d_2) = 0$, so that $x \cdot (b_1 \cdot d_1 + c_1 \cdot d_2) = c \cdot (b_1 \cdot d_1) + b \cdot (c_1 \cdot d_2)$. By (ii) $c \cdot (b_1 \cdot d_1), b \cdot (c_1 \cdot d_2) \in \zeta_1(B) \cap \zeta_1(C)$.

(v) By **Lemma 2.5** $L \bullet D = B \bullet D + C \bullet D$, so that

$L \bullet (L \bullet D) \leq L \bullet (B \bullet D + C \bullet D) \leq \zeta_1(L)$, therefore $L \bullet (L \bullet (L \bullet D) = \langle 0 \rangle$.

(vi) It follows from **Theorem 2.7.**

(vii) Futhermore, $d \cdot x = d \cdot (b + c) = d \cdot b + d \cdot c$. It follows that $D \bullet L \leq D \bullet B + D \bullet C$. Since $d \in B$, $d \cdot b \in B_1$. Since $d \in C$, $d \cdot c \in C_1$. On the other hand, $D$ is an ideal of $A$, and therefore, $d \cdot b \in B_1 \cap D$, $d \cdot c \in C_1 \cap D$, so that

$$D \bullet L \leq (D \cap \zeta_1(B)) + (D \cap \zeta_1(C)) \leq \zeta_1(D).$$

(viii) Let $d_1$ be an arbitrary element of $D \cap \zeta_1(B)$, $d_2$ be an arbitrary element of $D \cap \zeta_1(C)$. Since $d_1 \in D \leq C$, $d_1 \cdot d_2 = 0$, and similarly, $d_2 \cdot d_1 = 0$. Using these equalities, it is not hard to prove that $(D \cap B_1) + (D \cap C_1)$ is a subalgebra of $L$, and this subalgebra is abelian. Let $x = b + c$, $b \in B$, $c \in C$. We have

$$x \cdot (d_1 + d_2) = x \cdot d_1 + x \cdot d_2 = (b + c) \cdot d_1 + (b + c) \cdot d_2 =$$
$$b \cdot d_1 + c \cdot d_1 + b \cdot d_2 + c \cdot d_2 = c \cdot d_1 + b \cdot d_2.$$

Since $d_1 \in D \leq C$, $c \cdot d_1 \in C_1 \cap$ , and the fact that $D$ is an ideal implies that $c \cdot d_1 \in D$, so that $c \cdot d_1 \in C_1 \cap D$. Similarly, $b \cdot d_2 \in B_1 \cap D$. Thus, $x \cdot (d_1 + d_2) \in (D \cap \zeta_1(B)) + (D \cap \zeta_1(C))$.

Furthermore,

$$(d_1 + d_2) \cdot x = (d_1 + d_2) \cdot (b + c) = d_1 \cdot b + d_1 \cdot c + d_2 \cdot b + d_2 \cdot c = d_1 \cdot c + d_2 \cdot b.$$

Using the above arguments, we obtain that $d_1 \cdot c \in C_1 \cap D$, $d_2 \cdot b \in B_1 \cap D$, so that $(d_1 + d_2) \cdot x \in (D \cap \zeta_1(B)) + (D \cap \zeta_1(C))$. This means that $(D \cap \zeta_1(B)) + (D \cap \zeta_1(C))$ is an ideal of $L$. It is not hard to prove that this ideal is abelian.

(ix) Let $z$ be an arbitrary element of $\zeta(D)$ and $x$ be an arbitrary element of $A$. Then $x = b + c$ where $b \in B$, $c \in C$. We have $x \cdot z = b \cdot z + c \cdot z$. The inclusion $D \leq B$ and the fact that $D$ is an ideal imply that

$$b \cdot z \in B^2 \cap D \leq \zeta(B) \cap D \leq \zeta(D),\ c \cdot z \in C^2 \cap D \leq \zeta(C) \cap D \leq \zeta(D),$$

so that $x \cdot z \in \zeta(D)$.

We have $z \cdot x = z \cdot b + z \cdot c$. Using the above arguments, we obtain that $z \cdot b, z \cdot c \in \zeta(D)$. Thus, we obtain that $\zeta(D)$ is an ideal of $L$.

(x) Since $D$ is an ideal, $D \bullet L \leq D$. The inclusion $D \leq B$ implies that $D = \zeta_2(D)$. Let $u$ be an arbitrary element of $D \bullet L$. Then $u = \sum_{1 \leq j \leq n} d_j \cdot a_j$ for some elements $d_j \in D$, $a_j \in L$, $1 \leq j \leq n$. Moreover, $a_j = b_j + c_j$ for some elements $b_j \in B$, $c_j \in C$, $1 \leq j \leq n$, so that then

$$u = \sum_{1 \leq j \leq n} d_j \cdot a_j = \sum_{1 \leq j \leq n} d_j \cdot (b_j + c_j) =$$
$$\sum_{1 \leq j \leq n} (d_j \cdot b_j + d_j \cdot c_j) = \sum_{1 \leq j \leq n} d_j \cdot b_j + \sum_{1 \leq j \leq n} d_j \cdot c_j.$$

Let $x = b + c$ be an arbitrary element of $A$, $b \in B$, $c \in C$. We have

$$u \cdot x = u \cdot (b + c) = u \cdot (b + c) = u \cdot b + u \cdot c =$$
$$(\sum_{1 \leq j \leq n} d_j \cdot b_j + \sum_{1 \leq j \leq n} d_j \cdot c_j) \cdot b + (\sum_{1 \leq j \leq n} d_j \cdot b_j + \sum_{1 \leq j \leq n} d_j \cdot c_j) \cdot c =$$
$$\sum_{1 \leq j \leq n} ((d_j \cdot b_j) \cdot b) + (\sum_{1 \leq j \leq n} ((d_j \cdot c_j) \cdot b) +$$
$$\sum_{1 \leq j \leq n} ((d_j \cdot b_j) \cdot c) + (\sum_{1 \leq j \leq n} ((d_j \cdot c_j) \cdot c).$$

The inclusion $d_j \in D \leq B$ implies that $d_j \cdot b_j \in B^{(2)}$ and $(d_j \cdot b_j) \cdot b \in B^{(2)} \bullet B = B^{(3)} = \langle 0 \rangle$. Similarly, $(d_j \cdot c_j) \cdot c = 0$. Hence, we obtain that every element of $(D \bullet L) \bullet L$ is a sum of elements of one of the following forms: $(d \cdot c) \cdot b$, $(d \cdot b) \cdot c$, $b \in B$, $c \in C$, $d \in D$. In other words, $(D \bullet L) \bullet L \leq (D \bullet C) \bullet B + (D \bullet B) \bullet C$. Using the above arguments, we can prove the inclusions $(D \bullet C) \bullet B, (D \bullet B) \bullet C \leq \zeta(D)$. We have further

$$((d \cdot c_1) \cdot b_1) \cdot x = ((d \cdot c_1) \cdot b_1) \cdot (b + c) = ((d \cdot c_1) \cdot b_1) \cdot b + ((d \cdot c_1) \cdot b_1) \cdot c,$$
$$((d \cdot b_1) \cdot c_1) \cdot x = ((d \cdot b_1) \cdot c_1) \cdot (b + c) = ((d \cdot b_1) \cdot c_1) \cdot b + ((d \cdot b_1) \cdot c_1) \cdot c.$$

Since $D$ is an ideal of $L$, $d \cdot c_1 \in D \leq B$. Then $((d \cdot c_1) \cdot b_1) \in B^{(2)}$ and $((d \cdot c_1) \cdot b_1) \cdot b \in B^{(2)} \bullet B = B^{(3)} = \langle 0 \rangle$. Similarly, $((d \cdot b_1) \cdot c_1) \cdot c = 0$. Hence, we obtain that every element of $((D \bullet L) \bullet L) \bullet L$ is a sum of elements of one of the following forms:

$$((d \cdot c_1) \cdot b_1) \cdot c \in (D \bullet B) \bullet C,\ ((d \cdot b_1) \cdot c_1) \cdot b \in (D \bullet C) \bullet B,\ b_1, b \in B,\ c_1, c \in C,\ d \in D.$$

In other words, $((D \bullet L) \bullet L) \bullet L \leq (D \bullet B) \bullet C + (D \bullet C) \bullet B$. By what was proved above, $(D \bullet B) \bullet C + (D \bullet C) \bullet B \leq (D \cap \zeta_1(B)) + (D \cap \zeta_1(C)) \leq \zeta_1(D)$.

(xi) By **Corollary 1.5** the factor – algebra $L/D$ is nilpotent and **uzl**$(L/D) \leq 2$. It follows that $L^{(3)} \leq D$, and

$$L^{(4)} = L^{(3)} \bullet L \leq D \bullet L,\ L^{(5)} = L^{(4)} \bullet L \leq (D \bullet L) \bullet L \leq (D \bullet B) \bullet C + (D \bullet C) \bullet B \leq \zeta(D).$$

(xii) By (iii)

$$L \bullet D \leq (D \cap \zeta_1(B)) + (D \cap \zeta_1(C)) \leq \zeta_1(D)$$

and by (vii)

$$D \bullet L \leq (D \cap \zeta_1(B)) + (D \cap \zeta_1(C)) \leq \zeta_1(D).$$

This means that the factor $D/(D \cap \zeta_1(B)) + (D \cap \zeta_1(C)))$ is central. As we have noted above, $A/D$ is nilpotent and $\mathbf{uzl}(A/D) \leq 2$. It follows that the factor – algebra $L/(D \cap \zeta_1(B)) + (D \cap \zeta_1(C)))$ is nilpotent and

$$\mathbf{uzl}(L/(D \cap \zeta_1(B)) + (D \cap \zeta_1(C))) \leq 3.$$

Leonid.A. Kurdachenko, National Oles Honchar Dnipro University, Ukraine
lkurdachenko@gmail.com
Igor Ya. Subbotin, National University, USA
isubboti@nu.edu